\documentclass[11pt]{article}
\usepackage[a4paper, left=2cm, right=2cm, top=2cm, bottom=2cm]{geometry}
\usepackage{graphicx}
\usepackage{amssymb}
\usepackage{epstopdf}
\usepackage[hidelinks]{hyperref}
\title{New Ideas to the Design of Algorithms Based on  Derivatives}
\author{Fl\'avio Barbosa$^{\dagger,*}$ \& Fernando Nogueira$^{\dagger,**}$}
\date{$^{\dagger}$Universidade Federal de Juiz de Fora - Brazil}                                           

\begin{document}
\maketitle
\noindent $^{*}$ flavio.barbosa@ufjf.br - corresponding author  \\
\noindent $^{**}$fernando.nogueira@ufjf.br

\section*{Abstract}

This article proposes new perspectives for developing derivative-based numerical algorithms, supported by the introduction of Generalized Derivative operators. It demonstrates that these operators  have the potential to enhance and extend existing derivative-based numerical methods. To this end, two iterative derivative-driven methods are examined and refined: the Newton–Raphson method and the Gradient method. For both approaches, Generalized Derivatives are introduced with the goal of reducing the number of iterations required for convergence. These modifications are presented through geometric interpretations of the proposed constructions, which clearly illustrate their convergence-accelerating properties. The concluding remarks emphasize the significant opportunity to advance and refine numerical algorithms through the use of Generalized Derivatives.

\begin{center}
{\bf keywords:} Generalized Derivatives;
Derivative-Based Methods;
Gradient Method;
Newton–Raphson Method
\end{center}

\section{Introduction}
\label{introducao}

Derivative-based numerical methods constitute one of the most important class of computational tools in applied mathematics \cite{BurdenFaires}, optimization \cite{NocedalWright}, and scientific computing \cite{Quarteroni}. By exploiting local differential information—most commonly the first derivative—these methods generate iterative schemes capable of efficiently approximating roots, identifying extrema, and capturing structural behavior of nonlinear functions. Classical examples include Newton-type algorithms \cite{OrtegaRhynes}, gradient-based optimization strategies \cite{Bertsekas}, and various quasi-Newton and higher-order schemes \cite{DennisSchnabel}.

Within this broad spectrum of techniques, the Newton-Raphson (NR) method \cite{art9,art10,art11} and the Gradient (G) method \cite{Gagrani_2023, Biau_2019, Bottou_2018} stand out as two of the most influential and widely employed derivative-based algorithms. The NR method is distinguished by its fast local convergence and its elegant geometric interpretation based on tangent approximations. The G method, in turn, forms the backbone of continuous optimization, owing to its ability to leverage first-order information to guide descent along directions of steepest decrease. Both methods strike a compelling balance between theoretical rigor, computational feasibility, and practical performance, factors that underscore their enduring relevance across scientific \cite{Wibisono_2016, Bertsekas_2000} and engineering applications \cite{art8,anapaula}. Their ubiquity further motivates the search for extensions capable of enhancing their convergence properties while preserving their conceptual and computational simplicity.

To explore such extensions, this article applies the concept of a Generalized Derivative \cite{Nogueira_2024}, designed as a natural enlargement of Newton's classical notion of differentiation and aimed at broadening the mathematical tools available for algorithmic development.  

Based on this concept, it was developed Quadratic and Cubic Derivatives, derived  from the generalized differentiation framework. These new derivative forms provide enriched local information and enable the formulation of enhanced versions of the NR and G methods. The resulting algorithms are specifically designed to reduce the number of iterations required to achieve convergence. Their convergence-accelerating characteristics are justified through geometric constructions, which make explicit how the Generalized Derivatives reshape the underlying iterative dynamics.

The authors believe that introducing the concept of Generalized Derivatives has strong potential to enable the development and improvement of various numerical methods, particularly those that rely on derivative-based formulations.

This article is organized into five additional sections. Section 2 introduces the Generalized Derivative and presents the formulations of the Quadratic and Cubic Generalized Derivatives. Sections 3 and 4 then describe the proposed Quadratic Gradient and Cubic Newton–Raphson methods, respectively. Section 5 provides a numerical example that illustrates the performance gains achieved by these methods. Finally, Section 6 offers the final discussions and conclusions, summarizing the main findings and outlining promising directions for future research in numerical algorithms based on generalized differentiation.


\section{Generalized Derivative}


The Generalized Derivative \cite{Nogueira_2024}  retains a direct geometric interpretation while extending classical differentiation in a manner that naturally incorporates higher-order structures.

Let $y(x)$ be a real-valued function and let $\Delta \in \mathbb{R}$. Using the points $(x, y(x))$ and $(x+\Delta, y(x+\Delta))$, one can compute the angular coefficient $a_1$ and the linear coefficient $a_0$ of the linear equation $y(x) = a_1 x + a_0$ by solving the following linear system $S_1$:

\begin{equation} \label{sistema_linear1}
S_1:\left\{
\begin{array}{ll}
y(x)=a_1x+a_0,\\[0.2cm]
y(x+\Delta)=a_1(x+\Delta)+a_0,
\end{array}
\right.
\end{equation}

Considering that $a_0$ and $a_1$ depend on the chosen point $(x, y(x))$, the solution of system $S_1$ can be expressed as:
\begin{eqnarray}
a_0(x) & =  &  \frac{y(x)(x+\Delta)-y(x+\Delta)x}{\Delta}, \label{a_0} \\
a_1(x) & =  &  \frac{y(x+\Delta)-y(x)}{\Delta}. \label{a_1}
\end{eqnarray}

The tangent line to the graph of $y(x)$ can now be obtained by taking the limit in (\ref{a_0})–(\ref{a_1}) as $\Delta \to 0$:
\begin{eqnarray}
a_0(x) & =  &   \lim_{\Delta \to 0} \frac{y(x)(x+\Delta)-y(x+\Delta)x}{\Delta}, \label{a0a} \\
a_1(x) & =  &  \lim_{\Delta \to 0}   \frac{y(x+\Delta)-y(x)}{\Delta}. \label{a1a}
\end{eqnarray}

The variable $a_1(x)$ corresponds to the classical derivative and gives the angular coefficient of the tangent line to $y(x)$ at $x$. Similarly, $a_0(x)$ gives the linear coefficient of the tangent line, although this quantity is not explicitly used in classical differential calculus.

The concept of the Derivator Function ($\mathfrak{D}$) is now introduced. In this preliminary formulation, the Derivator Function serves as the foundational interpolation structure and provides a simplified means of computing the tangent approximation to $y(x)$. It is defined as:

\begin{equation}
\mathfrak{D} = a_1 x + a_0.
\end{equation}

The Derivator Function may be understood as the tangent approximation to $y(x)$ at which the Generalized Derivative is evaluated. In its simplest form, it corresponds to a linear function, yielding the Linear Derivative in direct analogy with the classical derivative. However, the construction is not restricted to linear forms; the Derivator Function may represent other classes of functions, such as a parabola or a higher-order interpolant. For example, by defining

\begin{equation}
\mathfrak{D} = a_2 x^2 + a_1 x + a_0,
\end{equation}

\noindent the resolution of system $S_2$:

\begin{equation} \label{sistema_par}
S_2:\left\{
\begin{array}{ll}
y(x)=a_2x^2+a_1x+a_0, \\[0.15cm]
y(x+\Delta)=a_2(x+\Delta)^2+a_1(x+\Delta)+a_0, \\[0.15cm]
y(x+2\Delta)=a_2(x+2\Delta)^2+a_1(x+2\Delta)+a_0,
\end{array}
\right.
\end{equation}

\noindent leads to coefficients $a_2$, $a_1$, and $a_0$ when $\Delta \to 0$:

\begin{eqnarray}
a_2(x) & =  &   \lim_{\Delta \to 0}\frac{1}{2} 
\frac{y(x)-2y(x+\Delta)+y(x+2\Delta)}{\Delta^2},  \label{a2q} \\
a_1(x) & =  &   \lim_{\Delta \to 0}\frac{1}{2} 
\frac{k_0y(x)-k_1y(x+\Delta)+k_2y(x+2\Delta)}{\Delta^2}, \label{a1q} \\
a_0(x) & =  &   \lim_{\Delta \to 0}\frac{1}{2} 
\frac{k_0y(x)+k_1y(x+\Delta)+k_2y(x+2\Delta)}{\Delta^2}, \label{a0q}
\end{eqnarray}

\noindent where $k_0= 2x+3\Delta$, $k_1= 4x+4\Delta$, and $k_2= 2x+\Delta$.

\begin{figure}[!h]
\centering
\includegraphics[width=9.0cm]{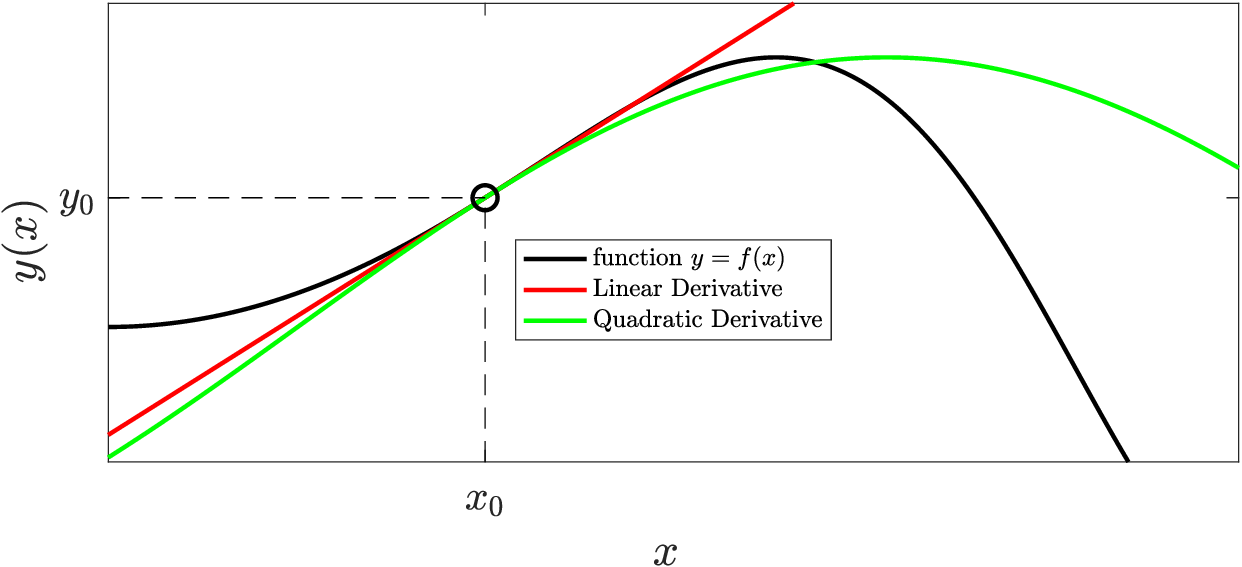}
\caption{Geometric interpretation of the Linear and Quadratic Derivatives}
\label{geometrica}
\end{figure}

Figure~\ref{geometrica} provides an intuitive illustration of the Generalized Derivative framework. The black curve depicts a generic function $y(x)$. At the point $(x_0,y_0)$, two distinct tangent constructions are shown: the red curve corresponds to the linear tangent associated with the classical derivative (Linear Derivative), while the green curve represents a parabolic tangent (Quadratic Derivative), capturing higher-order geometric information.

Following the same reasoning used to derive the Linear and Quadratic Derivatives, one may also construct the Cubic Derivative by choosing a cubic polynomial as the Derivator Function:

\begin{equation}
\mathfrak{D}=a_3 x^3 + a_2 x^2 + a_1 x + a_0.
\label{dersin}
\end{equation}

The resolution of system $S_3$:

\begin{equation} \label{sistema_C}
S_3:\left\{
\begin{array}{ll}
y(x)=y=a_3x^3+a_2x^2+a_1x+a_0,\\[0.15cm]
y(x+\Delta)=y_1=a_3(x+\Delta)^3+a_2(x+\Delta)^2+a_1(x+\Delta)+a_0,\\[0.15cm]
y(x+2\Delta)=y_2=a_3(x+2\Delta)^3+a_2(x+2\Delta)^2+a_1(x+2\Delta)+a_0,\\[0.15cm]
y(x+3\Delta)=y_3=a_3(x+3\Delta)^3+a_2(x+3\Delta)^2+a_1(x+3\Delta)+a_0,
\end{array}
\right.
\end{equation}

\noindent leads to coefficients $a_3$, $a_2$, $a_1$, and $a_0$ when $\Delta \to 0$:

\begin{equation}
a_3(x) =  \lim_{\Delta \to 0} \frac{3 y_1 - y - 3 y_2 + y_3}{6 \Delta^3}. 
\label{a3c}
\end{equation}

\begin{equation}
a_2(x) =  \lim_{\Delta \to 0} 
-\frac{5\Delta y_1 + 3x y_1 - 2\Delta y - 4\Delta y_2 + \Delta y_3 - x y - 3x y_2 + x y_3}{2\Delta^3}.
\end{equation}

\begin{equation}
a_1(x)= \lim_{\Delta \to 0} A_1,
\end{equation}

\noindent where (in compact form):

{\scriptsize
\[
A_1 = 
-\frac{
3 x^2 y + 9 x^2 y_2 - 3 x^2 y_3 - 18 \Delta^2 y_1 - 9 x^2 y_1 
+ 11 \Delta^2 y + 9 \Delta^2 y_2 - 2 \Delta^2 y_3 
- 30 \Delta x y_1 + 12 \Delta x y + 24 \Delta x y_2 - 6 \Delta x y_3
}{6 \Delta^3}.
\]
}

\begin{equation}
a_0(x)=\lim_{\Delta \to 0} A_0,
\label{a0c}
\end{equation}

\noindent with

{\scriptsize
\[
A_0 =
\frac{
x^3 y + 3 x^3 y_2 - x^3 y_3 - 3 x^3 y_1 + 6 \Delta^3 y 
+ 6 \Delta x^2 y + 11 \Delta^2 x y + 12 \Delta x^2 y_2 
+ 9 \Delta^2 x y_2 - 3 \Delta x^2 y_3 - 2 \Delta^2 x y_3 
- 15 \Delta x^2 y_1 - 18 \Delta^2 x y_1
}{6 \Delta^3}.
\]
}

The concept of a Generalized Derivative is therefore inherently dependent on the choice of the Derivator Function. Depending on this choice, the resulting Generalized Derivative may encode structural information that would not be explicitly accessible through classical differentiation.  
For instance, Nogueira and Barbosa~\cite{Nogueira_2024} demonstrate that for a chirp-type signal, the Cosenoidal Derivative ($\mathfrak{D} = \cos(\omega t + \phi)$) yields the instantaneous frequency $\omega$ and phase $\phi$.

In this article, the Quadratic and Cubic Derivatives are incorporated into the G and NR methods, respectively, with the goal of reducing the number of iterations required for convergence.

\section{Formulation of Quadratic G Method}

The G method is one of the most fundamental iterative procedures for locating local extrema of differentiable functions. Given a scalar function $y(x)$, the method constructs a sequence of iterates $x_n$ that progressively approach a local minimum or maximum by following the direction of steepest local variation of $y(x)$.

For a minimization problem, the update equation is
\begin{equation}
    x_{n+1} = x_n - a \, y^{'}_{n}, 
        \label{gradmin}
\end{equation}

\noindent where $a > 0$ is the step size and $y^{'}_{n} = \nabla y(x_n)$ denotes the gradient of the function evaluated at $x_n$. 

For a maximization problem, the direction of the gradient is reversed:
\begin{equation}
    x_{n+1} = x_n + a \, y^{'}_{n}.
    \label{gradmax}
\end{equation}

A central challenge in the practical application of the G method lies in the selection of an appropriate value for the parameter $a$. If $a$ is chosen too small, the iterative process converges very slowly, often requiring an impractically large number of iterations to show meaningful progress. Conversely, if $a$ is too large, the iterates may overshoot the region of interest, resulting in oscillations or even divergence. Thus, the efficiency and reliability of the G method  depend on a well-chosen step size.

To mitigate these issues, various variable step-size strategies have been developed. Among the most widely used are line-search procedures, such as Armijo backtracking\cite{Armijo1966}, and other adaptive mechanisms \cite{Davis_2025,Krutikov_2024} that adjust the step size at each iteration to balance convergence speed and stability \cite{NocedalWright}.

In this article, a  variant of the classical G method is presented, aiming to significantly reduce 
the number of iterations required to converge to a local extremum of the function $ y(x) $. 
The proposed approach replaces the standard Linear Derivatives, referred to in this text as the L-G method, 
with Quadratic Derivatives. This modification yields the Quadratic G (Q-G) method and increases the method’s 
sensitivity to local curvature, thereby improving convergence efficiency.

\begin{figure}[h]
\centering
       \includegraphics[width=12.0cm] {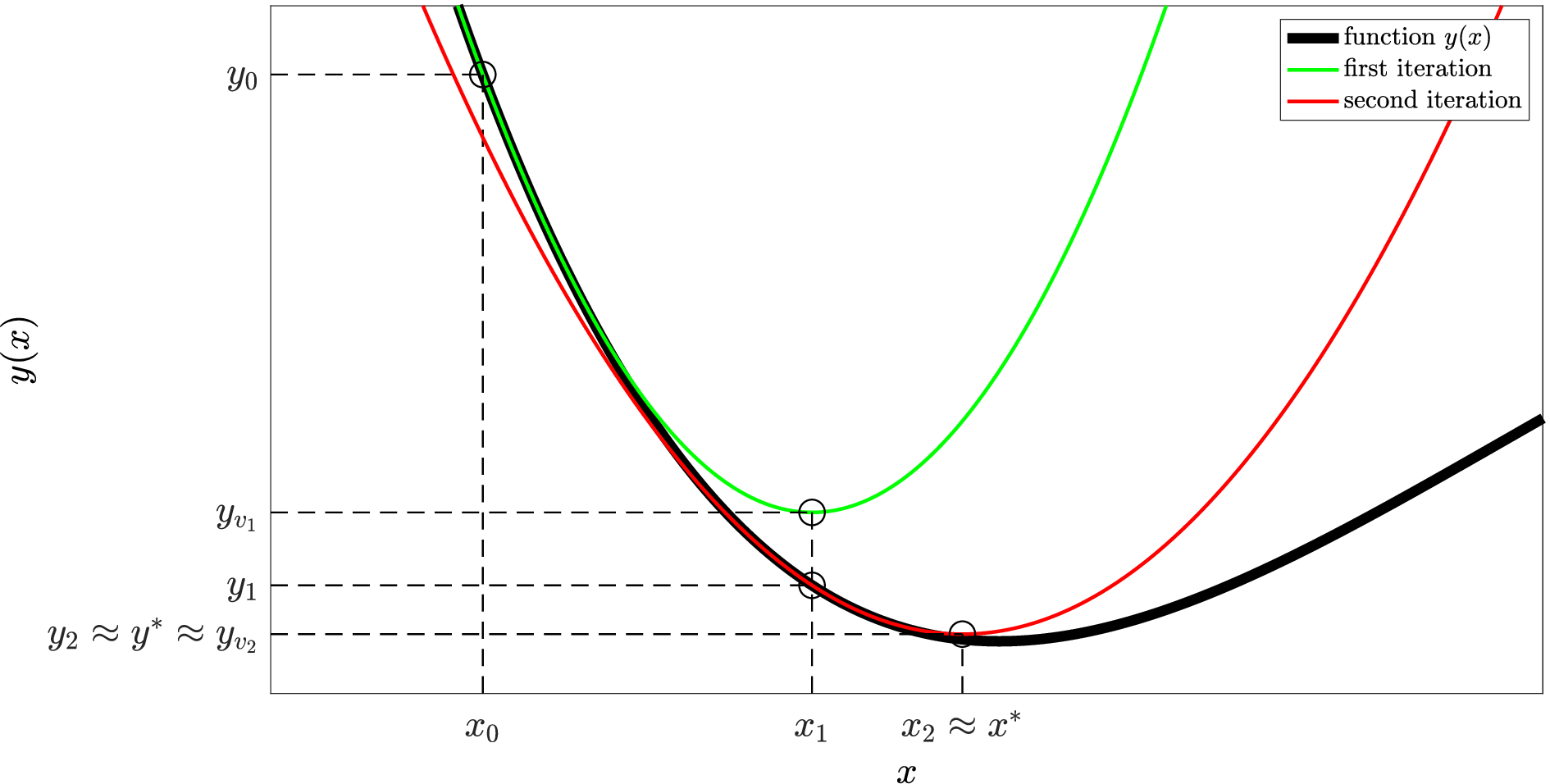} 
       \caption{Graphical interpretation of Q-G method towards $({x}^*,y^*)$}
       \label{G}
\end{figure}

The geometric scheme of the proposed modification is shown in Figure~\ref{G}. Starting from an initial estimate $(x_0, y_0)$, the Quadratic Derivative is defined through this point (green curve), establishing the first iteration toward the local minimum $(x^*, y^*)$. The vertex of the Quadratic Derivative, that is, the point $(x_1, y_{v_1})$ on the green curve, constitutes the second estimate of $(x^*, y^*)$. For the second iteration, the procedure starts from the point $(x_1, y_1)$, and a new Quadratic Derivative is defined through this point (red curve). The vertex of the red parabola, namely $(x_2, y_{v_2})$, provides the next estimate of $(x^*, y^*)$. This iterative process can be repeated until a stopping criterion is met.  Graphically, it can be observed that each iteration yields a closer approximation to $(x^*, y^*)$. 

The update equation for the Q-G method is:

\begin{equation}
    x_{n+1} =-\frac{a_1(x_n)}{a_2(x_n)},
\end{equation}

\noindent and
\begin{equation}
    y_{n+1} =-\frac{[a_1(x_n)]^2 - 4 a_2(x_n) a_0(x_n)}{4 a_2(x_n)},
\end{equation}

\noindent where $a_2(x_n)$, $a_1(x_n)$ and $a_0(x_n)$ are given by Equations~\ref{a2q}, \ref{a1q} and \ref{a0q}, respectively.

Other types of Derivator Functions may be employed, and identifying the most suitable one for each specific case remains one of the many open directions for future research on the topic. When prior knowledge about the function \( y(x) \), whose local extrema are being sought, is available, selecting a Derivator Function that closely resembles \( y(x) \) can significantly accelerate convergence.

\section{Formulation of Cubic NR Method}

The NR method approximates the roots of a given function through iterative refinement of an initial estimate, exploiting the local behavior of the function via its Linear Derivative to drive the sequence toward convergence. In practical terms, the method computes the next iterate $x_{n+1}$ as the root of the tangent line to the function $y(x)$ at the point $(x_n, y_n)$, repeating this procedure until convergence is achieved.

Since the tangent curve used in the classical NR method is a straight line, obtained through a Linear Derivator Function, this article refers to the classical formulation as the Linear NR (L-NR) method.

Once the tangent curve to the function $y(x)$ passing through the point $(x_0, y_0)$, the initial guess, is determined, the next step is to compute its root. In the case of the L-NR method, this calculation is straightforward, leading to the iterative equation:

\begin{equation}
x_{n+1}=x_{n}-\frac{y_{n}}{y^{'}_{n}}.
\end{equation}

Instead of relying on a Linear Derivator Function, the authors propose the use of a Cubic Derivator Function to guide the initial estimate toward the root. In this formulation, determining the roots of the tangent curve becomes considerably more complex than in the linear case. Nevertheless, by applying Cardano's closed-form solution \cite{cardano} (see Appendix), the three roots of the cubic polynomial can be computed, leading to the iterative equation that defines the Cubic NR (C-NR) method:

\begin{equation}
x_{n+1} = \mathrm{RootOf}\!\left[a_3(x_{n})\,x^{3} + a_2(x_{n})\,x^{2} + a_1(x_{n})\,x + a_0(x_{n})\right],
\label{cubica}
\end{equation}

\noindent where $\mathrm{RootOf}[\cdot]$ denotes the real root closest to the current estimate $x_{n}$, as recommended by the authors.



\begin{figure}[h]
\centering
       \includegraphics[width=12.0cm] {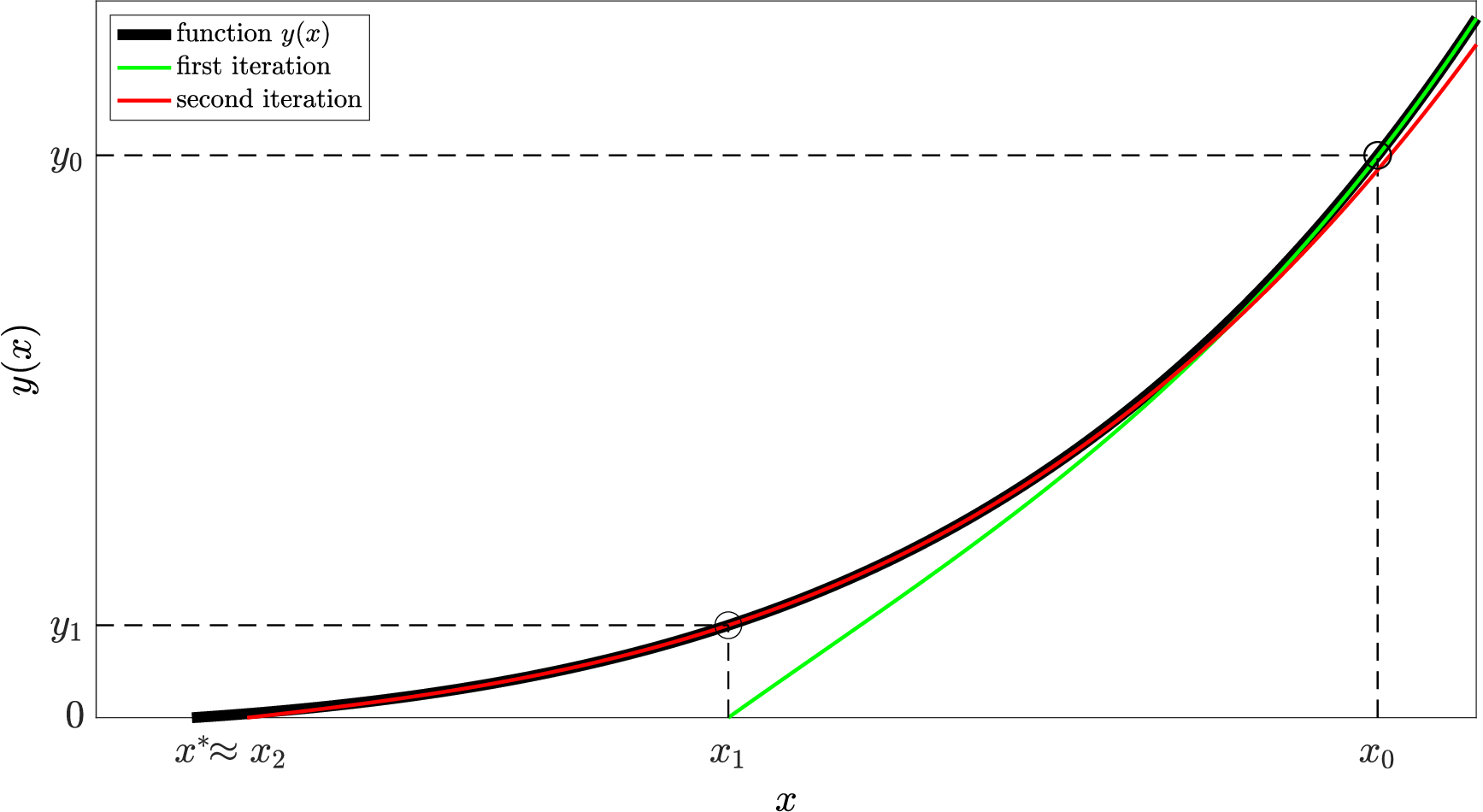} 
       \caption{Graphical interpretation  of C-NR method towards the root $x^*$}
       \label{NR}
\end{figure}

Figure~\ref{NR} presents a graphical interpretation of the C-NR method for the first two iterations in the search for the root $x^{*}$ of the function $y(x)$. Starting from the initial guess $(x_{0}, y_{0})$, the C-NR method yields the update $x_{1}$. The second iteration is then obtained by evaluating the Cubic Derivator Function at the point $(x_{1}, y_{1})$, resulting in the new approximation $x_{2}$. Notably, after only two iterations, the method already provides a remarkably approximation of the actual root $x^{*}$.

At this point of the text, it becomes possible to describe the Generalized NR method. In this scenario, the Derivator Function used for tangent calculation can be others, not solely a linear or cubic function, as in the case of the L-NR or C-NR methods, presented in this article.  Naturally, the selection of the Derivator Function will directly influence the convergence process, potentially speeding it up, slowing it down, or even rendering it impractical.

Similarly to the previous case, when there is prior knowledge of the function whose roots are sought, choosing a Derivator Function that closely matches it can significantly speed up convergence. In a limit situation, when the target function is cubic, using a cubic Derivator Function may lead to convergence in a single iteration, since its roots can be obtained directly from the closed-form solution of the cubic equation.

In all cases, it is advisable that the chosen Derivator Function possess at least one real root. When using a Quadratic Derivator Function, real roots may not always exist (see Eqs.~\ref{a2q}--\ref{a0q}). In such situations, the Quadratic NR method may converge if real roots are present; however, if the Quadratic Derivator Function has no real roots, the iterative process will necessarily fail to converge.

\section{Example of application}

The purpose of this application example is to demonstrate that employing the C-NR method for root finding and the Q-G method for determining local extrema leads to a substantially reduced number of iterations required for convergence. To this end,  the function:

\begin{equation}
y(x)=x^4 - 21x^3 + 149x^2 - 419x + 290
\label{func}
\end{equation}
   
\noindent  is used. This function exhibits three local extrema: local maxima at \( x = 2.60 \) and \( x = 8.28 \), 
and a local minimum at \( x = 4.87 \). It also has real roots at \( x = 1 \) and \( x = 10 \).

For the evaluations, analytical functions defining the Linear, Quadratic, and Cubic Derivatives were used. 
In all analyses, the convergence criterion adopted to determine convergence to a point of interest \( x^* \) 
at the \( n \)-th iteration was

\begin{equation}
|x_n - x_{n-1}| \leq 10^{-4},
\end{equation}

\noindent which leads to

\begin{equation}
x^* \approx x_n.
\end{equation}

\subsection{Evaluation of Q-G Method}

As previously noted, the choice of the parameter $a$ used in the L-G method is a key aspect of its performance.  Classical references recommend values of $a$ ranging from $10^{-3}$ to $10^{-1}$ \cite{NocedalWright}. In this work, the value $a = 0.05$ was adopted, while it is worth noting that the Q-G method does not require any parameter to be tuned.

Additionally, since the application of the L-G method requires specifying the type of extremum being investigated, Equation \ref{gradmin} is adopted here to target the minimization of $y(x)$. In contrast, such information is not required by the Q-G method, since it is capable of inherently identifying all extrema of the function. Naturally, through the use of appropriate constraints, one may still restrict the procedure to a desired form of extremization—either maximization or minimization.

\begin{figure}[h]
\centering
\includegraphics[width=0.8\textwidth]{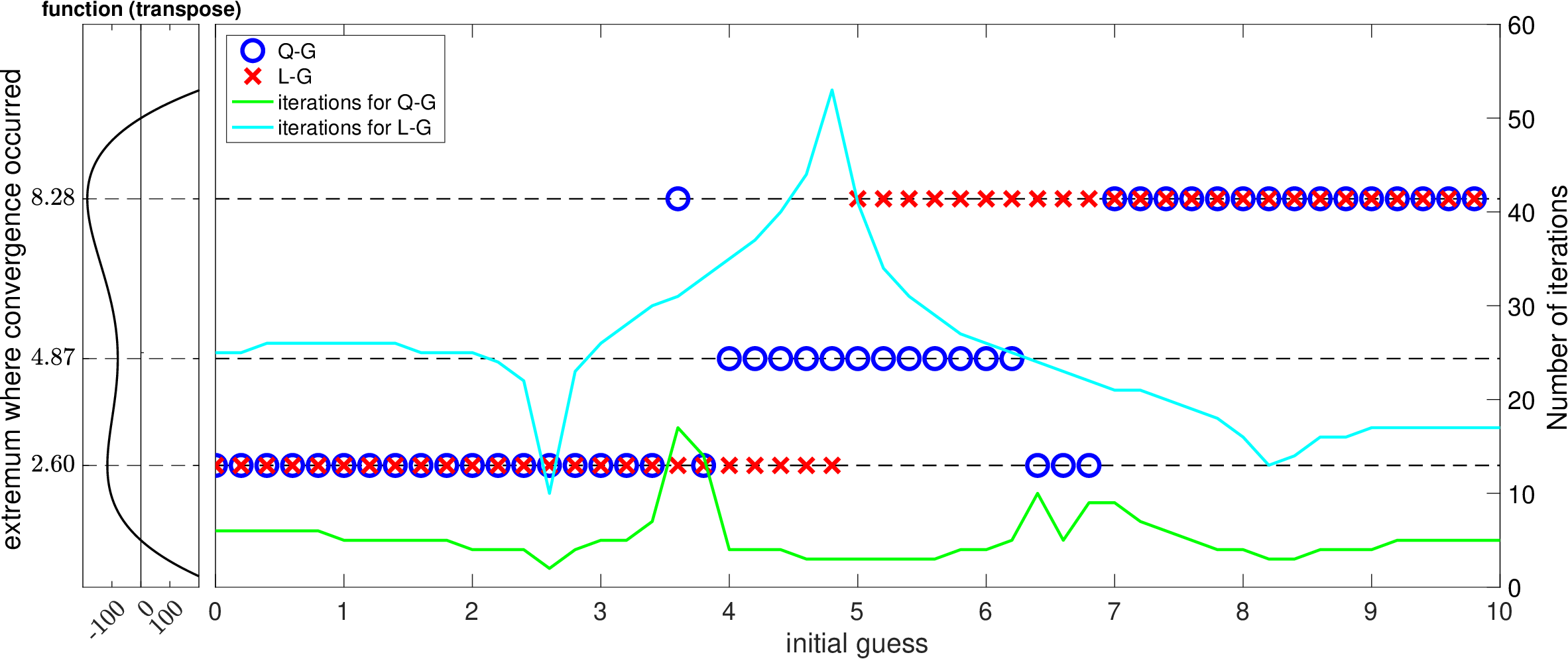}
\caption{Comparisons for function $y_c(x)$}
\label{f3}
\end{figure}

Figure~\ref{f3} presents a comparison between the results obtained using the Q-G and L-G methods. 
In this Figure, the left panel shows the graph of the function \( y(x) \) (Equation~\ref{func}), 
displayed in a transposed layout in which the vertical axis corresponds to the \( x \)-axis and 
the horizontal axis corresponds to the \( y \)-axis. The dashed lines indicate the extrema of 
the function \( y(x) \).  The main plot illustrates the initial guesses used for the computation of extrema along the 
horizontal axis and the corresponding extremum values obtained along the left vertical axis. 
Blue circles and red crosses denote the convergence points produced by the Q-G and L-G methods, 
respectively.  On the right side of the figure, the secondary vertical axis represents the number of iterations 
required to achieve convergence, shown using the green (Q-G) and cyan (L-G) series.

From Figure~\ref{f3}, it can be observed that both the Q-G and L-G algorithms converged to local extrema of $y(x)$. The L-G method identified the two local minima, whereas the Q-G method successfully identified all local extrema of $y(x)$. However, convergence to each extremum consistently required a substantially smaller number of iterations for the Q-G method across all initial estimates. On average, the L-G method required approximately five times more iterations than Q-G (25 versus 5.21), and the maximum number of iterations reached 53 and 17, respectively. Although $y(x)$ contains relatively flat regions that typically hinder iterative schemes, these features did not introduce significant difficulties for the Q-G method, in contrast to the behavior observed with the L-G method.

In addition, the L-G method exhibited a standard deviation in the number of iterations that was nearly three times larger than that of the Q-G method (8.29 versus 2.64). This increased variability indicates that L-G is more sensitive to the choice of initial conditions than Q-G.  Consequently, the Q-G method demonstrated greater robustness, as its performance remained more uniform across different starting points. The low variability observed for Q-G further facilitates the prediction of computational cost, making execution time and resource requirements easier to estimate in advance. It also reduces the likelihood of divergence or excessively slow convergence across different initial guesses. 

\subsection{Evaluation of C-NR Method}

\begin{figure}[h]
\centering
\includegraphics[width=0.8\textwidth]{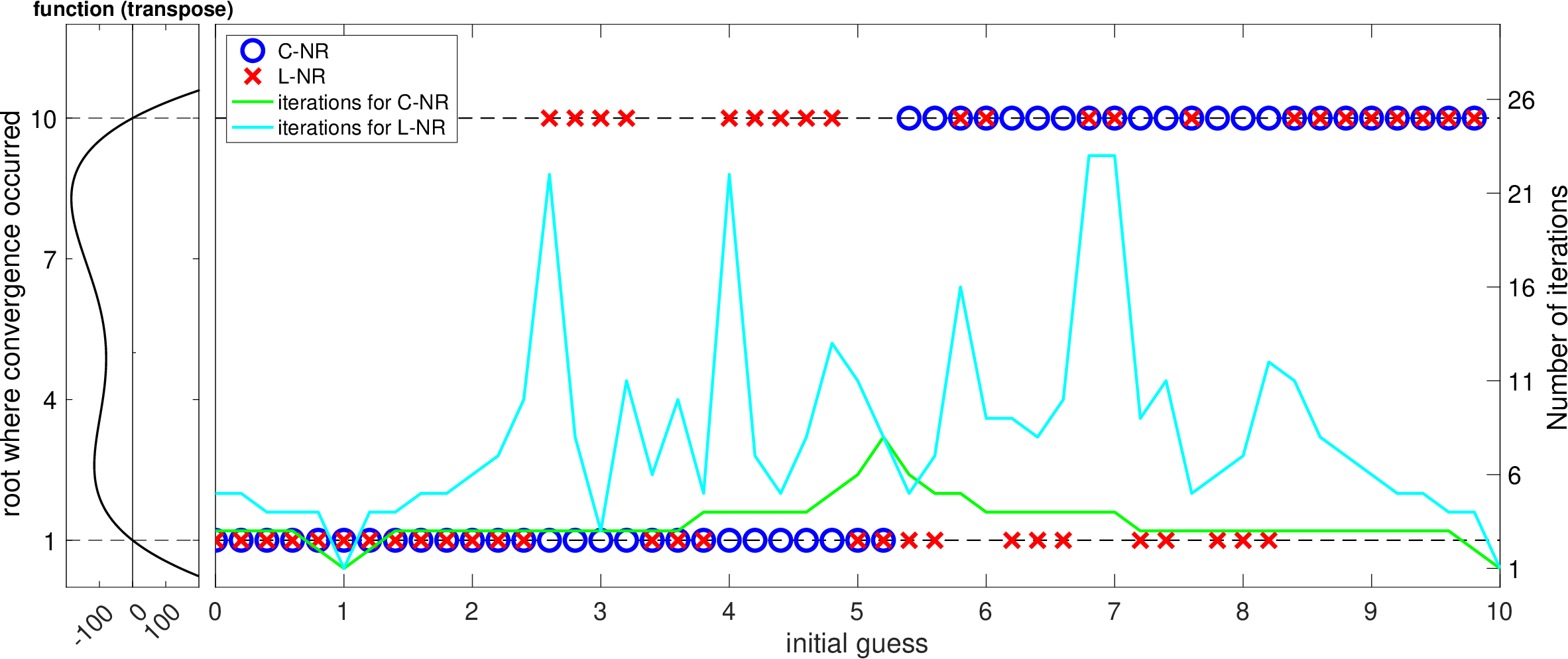}
\caption{Comparisons for function $y_b(x)$}
\label{f2}
\end{figure}

Figure~\ref{f2} presents a comparison between the results obtained using the C-NR and L-NR methods. Both methods converged to one of the roots of the problem for all initial values evaluated. However, the C-NR method consistently converged to the root closest to the initial estimate, without oscillating or switching between different roots, which indicates greater local stability and predictability. This behavior suggests that C-NR possesses better-defined basins of attraction and a reliable update mechanism that is not overly sensitive to perturbations in the initial guess.

From this figure, it can be observed that the C-NR method required a significantly smaller number of iterations than the L-NR method to achieve convergence for almost all initial estimates. On average, the L-NR method required 8.11 iterations, whereas the C-NR method required 3.41, with peak values of 23 and 8 iterations, respectively.  This result indicates that the C-NR method achieves a higher convergence rate compared to the L-NR method, allowing it to reach the solution in fewer iterations. 

The standard deviations observed for the number of iterations required to achieve convergence (1.17 for C-NR and 5.18 for L-NR) lead to remarks similar to those previously identified in the comparison between the Q-G and L-G methods. The  lower variability of the C-NR method indicates reduced sensitivity to the initial guess, and its performance therefore remains more consistent across different starting points. This behavior reflects a higher degree of robustness compared to the L-NR method.

\section{Final Discussions and Conclusions}

This work introduced a generalized framework for derivative operators and examined its implications for the development of numerical algorithms. Rather than focusing on validating new techniques through extensive test batteries, the main objective was to open a conceptual avenue and show that Generalized Derivatives offer a flexible and still largely unexplored structure capable of refining and accelerating classical iterative schemes. By revisiting two well-established methods, Newton–Raphson and Gradient-based approaches, the study demonstrated that incorporating Generalized Derivatives can significantly reduce the number of iterations needed for convergence while also improving robustness with respect to initial conditions.

The results obtained for the Q-G and C-NR methods reinforce this perspective. In both cases, the generalized formulations consistently outperformed their classical counterparts, achieving faster convergence even in regions where traditional methods tend to slow down or become less stable. They also exhibited lower sensitivity to initial guesses, lsuggesting better-defined attraction basins and a more predictable update mechanism.  Together, these observations show that Generalized Derivatives can enhance both efficiency and stability in iterative numerical schemes and point to their potential for enabling more reliable and effective algorithms.

 However, the principal purpose of this article extends beyond proposing two refined algorithms. The broader intention is to motivate further research by illustrating the rich landscape of opportunities that Generalized Derivatives open for numerical analysis.

Several promising directions naturally emerge:
\begin{itemize}
    \item \textbf{Exploration of diverse Derivator Functions.} Different classes of Generalized Derivatives may be better suited to specific types of optimization or root-finding problems. Systematic analyses of these operator families may reveal problem-dependent optimal choices.
    
    \item \textbf{Extension to multivariate settings.} The development of Generalized-Derivative-based algorithms for systems of equations and multivariate optimization represents a fertile research area, potentially enabling more efficient and robust high-dimensional solvers.
    
    \item \textbf{Generalized derivative-based series expansions.} Constructing new series analogous to Taylor expansions, but grounded in Generalized Derivatives, may offer improved local approximations and new analytical tools for function representation.
    
    \item \textbf{Rigorous convergence analysis.} A comprehensive theoretical study of the convergence properties of Generalized-Derivative-based algorithms is essential for establishing guarantees, understanding stability regimes, and guiding the construction of future methods.
    \item \textbf{Non-euclidean geometry analysis.} Employing elliptic or hyperbolic Derivator Functions may facilitate both the formulation and the interpretation of elliptic and hyperbolic geometries. This approach could yield a more direct framework than the classical derivative, whose structure is fundamentally rooted in Euclidean geometry and therefore less suited for describing non-Euclidean spaces.
    
\end{itemize}

\vspace{0.3cm}
\noindent {\bf Acknolegements} 
\vspace{0.1cm}

\noindent The authors would like to thank CNPq  (Conselho Nacional de Desenvolvimento Cient\'{\i}fico e Tecnol\'ogico) - grants 407256/2022-9 and  303550/2025-2 -  for their financial support.

\vspace{0.3cm}

\noindent {\bf Declaration of interest statement} 
\vspace{0.1cm}

\noindent The authors declare that they have no conflicts of interest related to this work.

\section*{Appendix: Summary of Cardano's Equations}

The classical method for solving a general cubic equation was developed by 
Gerolamo Cardano in the 16\textsuperscript{th} century. 
Consider the general cubic polynomial

\begin{equation}
ax^{3} + bx^{2} + cx + d = 0, \qquad a \neq 0.
\label{cubic_general}
\end{equation}

The first step consists of reducing \ref{cubic_general} to a 
depressed cubic, i.e., a cubic equation without the quadratic term.
By performing the substitution

\[
x = t - \frac{b}{3a},
\]

\noindent one obtains the depressed cubic of the form

\begin{equation}
t^{3} + pt + q = 0,
\label{depressed_cubic}
\end{equation}

\noindent  where the coefficients $p$ and $q$ are given by

\begin{eqnarray}
p &= &\frac{3ac - b^{2}}{3a^{2}}, \\
q &= & \frac{2b^{3} - 9abc + 27a^{2}d}{27a^{3}}.
\end{eqnarray}

The discriminant of the cubic is defined as

\begin{equation}
\Delta = \left(\frac{q}{2}\right)^{2} 
       + \left(\frac{p}{3}\right)^{3}.
\end{equation}

Cardano’s closed-form solution for the depressed cubic 
\ref{depressed_cubic} is

\begin{equation}
t = \sqrt[3]{-\frac{q}{2} + \sqrt{\Delta}}
  + \sqrt[3]{-\frac{q}{2} - \sqrt{\Delta}}.
\label{eq:cardano_t}
\end{equation}

Finally, the solution of the original cubic equation \ref{cubic_general}
follows from the back-substitution

\begin{equation}
x = t - \frac{b}{3a}.
\end{equation}

When $\Delta > 0$, equation \ref{eq:cardano_t} provides one real root 
and two complex conjugate roots. When $\Delta = 0$, all roots are real and 
at least two coincide. For $\Delta < 0$, all three roots are real, and 
trigonometric forms of Cardano’s formula are typically employed to avoid 
complex intermediate computations.

This set of expressions, collectively known as \textit{Cardano's equations}, 
remains one of the most celebrated closed-form solutions in algebra.

\end{document}